\documentclass{article}
\usepackage{amsxtra,amssymb,amsthm,amsmath,latexsym}

\theoremstyle{plain}

\newtheorem{theorem}{Theorem}[section]

\numberwithin{equation}{section}

\def\R{{\mathbb R}}

\def\dotu{\dot{u}}

\def\oH{{\overset{\circ}{H}}}
\def\oH1{{\overset{\circ}{H}\kern-.02in{}^1}}

\def\bee{\begin{equation*}}
\def\eee{\end{equation*}}
\def\be{\begin{equation}}
\def\ee{\end{equation}}
\newcommand{\pn}{\par\noindent}

\title{Implicit Function Theorem via the DSM}

\author{A G Ramm\\
\small Department of Mathematics\\[-0.8ex]
\small Kansas State University, Manhattan, KS 66506-2602, USA\\[-0.8ex]
\small \texttt{ramm@math.ksu.edu}\\
}

\begin{document}
\date{}
\maketitle

\begin{abstract}
Sufficient conditions are given for an  implicit function theorem
to hold. The result is established by an application of the
Dynamical Systems Method (DSM). It allows one to solve a class of
nonlinear operator equations in the case when the Fr\'echet
derivative of the nonlinear operator is a smoothing operator, so
that its inverse is an unbounded operator.
\end{abstract}
\pn{
 }


\section{Introduction}\label{S:1}
The aim of this paper is to demonstrate the power of the Dynamical
Systems Method (DSM) as a tool for proving theoretical results. The
DSM was systematically developed in \cite{R499} and applied to
solving nonlinear operator equations in \cite{R499} (see also
\cite{R485}), where the emphasis was on convergence and stability of
the DSM-based algorithms for solving operator equations, especially
nonlinear and ill-posed equations. The DSM for solving an operator
equation $F(u)=h$ consists of finding a nonlinear map $u\mapsto
\Phi(t,u)$, depending on a parameter $t\in [0,\infty)$, that has the
following three properties:

(1) the Cauchy problem
$$\dot u=\Phi(t,u), \quad u(0)=u_0
\quad (\dotu:=\frac{du(t)}{dt})$$
has a unique global solution $u(t)$ for a
given initial approximation
$u_0$;

(2) the limit $u(\infty)=\lim_{t\to
\infty}u(t)$ exists; and

(3) this limit solves the original equation $F(u)=h$, i.e.,
 $F(u(\infty))=h.$

\noindent The operator $F:H\to H$ is a nonlinear map in a
Hilbert space $H$.
It is assumed that the equation $F(u)=h$ has a solution, possibly
nonunique.

The problem is to find a $\Phi$ such that the
properties
(1), (2), and (3) hold.
Various choices of $\Phi$ for which these properties hold
are proposed in \cite{R499}, where  the DSM is
justified for wide classes of operator equations, in particular, for
some
classes of nonlinear ill-posed equations (i.e., equations $F(u)=0$
for which
the linear operator $F'(u)$ is not boundedly invertible). By
$F'(u)$ we denote the Fr\'echet derivative of the nonlinear map $F$
at the element $u$.

 In this note the DSM is used as a tool for proving a  "hard"
implicit function theorem.

Let us first recall
the usual implicit function theorem. Let $U$ solve the equation 
$F(U)=f$.

{\it Proposition:} {\it If $F(U)=f$, $F$ is a
$C^1$-map in a Hilbert space $H$, and $F'(U)$ is a boundedly
invertible operator, i.e., $\|[F'(U)]^{-1}]\|\leq m,$ then the equation
\be\label{e1} F(u)=h \ee
is uniquely solvable for every $h$
sufficiently close to $f$.}

For convenience of the reader we
include a  proof of this known result.

{\it Proof of the Proposition.} First, one can reduce the problem to the
case
$u=0$
and $h=0$.
This is done as follows. Let $u=U+z$,
$h-f=p$, $F(U+z)-F(U):=\phi(z).$ Then $\phi(0)=0$,
$\phi'(0)=F'(U)$, and equation \eqref{e1} is equivalent to the
equation \be\label{e2} \phi(z)=p, \ee with the assumptions
\be\label{e3} \phi(0)=0,\quad \lim_{z\to
0}\|\phi'(z)-\phi'(0)\|=0,\quad \|[\phi'(0)]^{-1}\|\leq m. \ee We
want to prove that equation \eqref{e2} under the assumptions
\eqref{e3} has a unique solution $z=z(p)$, such that $z(0)=0$, and
$\lim_{p\to 0}z(p)=0.$ To prove this, consider the equation
\be\label{e4} z=z-[\phi'(0)]^{-1}(\phi(z)-p):=B(z), \ee and check
that the operator $B$ is a contraction in a ball
$\mathcal{B}_\epsilon:=\{z:\ \|z\|\leq \epsilon\}$ if $\epsilon>0$
is sufficiently small, and $B$ maps $\mathcal{B}_\epsilon$ into
itself. If this is proved, then the desired result follows from the
contraction mapping principle.

One has \be\label{e5}
\|B(z)\|=\|z-[\phi'(0)]^{-1}(\phi'(0)z+\eta -p)\|\leq m\|\eta\|+m\|p\|,
\ee
where $\|\eta\|=o(\|z\|)$.
If $\epsilon$ is so small that
$m\|\eta\|<\frac{\epsilon}{2}$ and $p$ is so small that
$m\|p\|<\frac{\epsilon}{2}$, then $\|B(z)\|<\epsilon,$ so
$B:\mathcal{B}_\epsilon\to \mathcal{B}_\epsilon$.

Let us check that $B$ is a  contraction mapping in $\mathcal{B}_\epsilon$.
One has:
\be\label{e6}\begin{split}
\|Bz-By\|&=\|z-y-[\phi'(0)]^{-1}(\phi(z)-\phi(y))\|\\
&=\|z-y-[\phi'(0)]^{-1}\int_0^1\phi'(y+t(z-y))dt(z-y)\|\\
&\leq m\int_0^1\|\phi'(y+t(z-y))-\phi'(0)\|dt\|z-y\|.\end{split}\ee
If $y,z\in \mathcal{B}_\epsilon$, then
$$\sup_{0\leq t\leq
1}\|\phi'(y+t(z-y))-\phi'(0)\|=o(1), \qquad
\epsilon\to 0.$$
Therefore, if $\epsilon$ is so small that
$m o(1)<1$, then $B$ is a contraction mapping in
$\mathcal{B}_\epsilon$, and equation \eqref{e2} has a unique
solution $z=z(p)$ in $\mathcal{B}_\epsilon$, such that $z(0)=0.$
The proof is complete. \hfill $\Box$

The crucial assumptions, on which this proof is based, are assumptions
\eqref{e3}.

Suppose now that
$\phi'(0)$ {\it is not boundedly invertible, so that the last assumption
in \eqref{e3} is not valid}. Then a theorem which still guarantees
the existence of a solution to equation \eqref{e2} for some set of
$p$ is called a "hard" implicit function theorem. Examples of such
theorems one may find, e.g.,  in \cite{AG}, \cite{B}, \cite{D}, and
\cite{H}.

Our goal in this paper is to
establish a new theorem of this type using a new method of proof,
based on the Dynamical Systems Method (DSM). In \cite{R489}
we have demonstrated a theoretical application of the DSM by
establishing some surjectivity results for nonlinear 
operators.

The result, presented in this paper, is a new illustration of
the applicability of the DSM as a tool for proving theoretical results.

To formulate the
result, let us introduce the notion of a scale of Hilbert 
spaces $H_a$
(see \cite{K}). Let $H_a\subset H_b$ and $\|u\|_b\leq \|u\|_a$ if
$a\geq b$. Example of spaces $H_a$ is the scale of Sobolev spaces
$H_a=W^{a,2}(D)$, where $D\subset\R^n$ is a bounded domain with a
sufficiently smooth boundary.

Consider equation \eqref{e1}. Assume
that
\be\label{e7} F(U)=f;\qquad
F:H_a\to H_{a+\delta}, \qquad u\in B(U,R):=B_a(U,R), \ee
where $B_a(U,R):=\{u\ :\ \|u-U\|_a\leq R\}$ and $ \delta=const>0$,
and the operator $F: H_a\to H_{a+\delta}$ is continuous.
Furthermore, assume that
$A:=A(u):=F'(u)$ exists and is an isomorphism of $H_a$ onto  
$H_{a+\delta}$:
\be\label{e8} c_0\|v\|_a\leq \|A(u)v\|_{a+\delta}\leq
c'_0\|v\|_a,\qquad u,v\in B(U,R), \ee
that
\be\label{e9}\|A^{-1}(v)A(w)\|_a\leq c,\qquad v,w\in B(U,R),\ee
and
\be\label{e10} \|A^{-1}(u)[A(u)-A(v)]\|_a\leq c||u-v||_a,\qquad u,v\in 
B(U,R).
\ee
Here and below we denote by $c>0$ various constants. Note that
\eqref{e8} implies 
$$||A^{-1}(u)\psi||_a\leq c_0^{-1}||\psi||_{a+\delta},
\qquad \psi= A(u)[F(v)-h],\quad v\in B(U,R).$$
Assumption \eqref{e8} implies that $A(u)$ is a smoothing 
operator similar 
to a smoothing integral operator, and its inverse is similar to
the differentiation operator of order $\delta>0$. 
Therefore, the operator
$A^{-1}(u)=[F'(u)]^{-1}$ causes the "loss of the derivatives".
In general, this may lead to a breakdown of the Newton process
(method) in a finitely many steps. 
Our assumptions \eqref{e7}-\eqref{e10}  
guarantee that this will not happen.

Assume that
\be\label{e11}u_0\in B_a(U,\rho), 
\qquad h\in B_{a+\delta}(f,\rho),\ee 
where
$\rho>0$ is a sufficiently small number: 
$$\rho\leq\rho_0:= \frac{R}{1+c_0^{-1}(1+c_0')},$$ 
and $c_0, c_0'$ are the constants from \eqref{e8}.
Then $F(u_0)\in  B_{a+\delta}(f, c_0'\rho)$, because
$||F(u_0)-F(U)||\leq c_0'||u_0-U||\leq c_0'\rho$.

Consider the problem \be\label{e13}
\dot{u}=-[F'(u)]^{-1}(F(u)-h),\quad u(0)=u_0. \ee
Our basic result is:
\begin{theorem}\label{T:1}If the assumptions 
\eqref{e7}-\eqref{e11}
hold, and 
$0<\rho\leq\rho_0:=\frac{R}{1+c_0^{-1}(1+c_0')},$ 
where $c_0, c_0'$ are the
constants from \eqref{e8}, then problem \eqref{e13} has a 
unique
global solution $u(t)$, there exists $V:=u(\infty)$, \be\label{e14}
\lim_{t\to \infty}\|u(t)-V\|_a=0,\ee and \be\label{e15} F(V)=h.\ee
\end{theorem}
Theorem \ref{T:1} says that if $F(U)=f$ and $\rho\leq\rho_0$, then for any
$h\in B_{a+\delta}(f,\rho)$ equation \eqref{e1} is solvable and a solution
to \eqref{e1} is $u(\infty)$, where $u(\infty)$ solves problem
\eqref{e13}.

In Section \ref{S:2} we prove Theorem \ref{T:1}.

\section{Proof}\label{S:2}
Let us outline the ideas of the proof. The local existence and
uniqueness of the solution to \eqref{e13} will be established if one
verifies that the operator $A^{-1}(u)[F(u)-h]$ is 
locally Lipschitz in 
$H_a$. The global existence of this solution $u(t)$ will be 
established if one proves the uniform boundedness of $u(t)$:
\be\label{e16} \sup_{t\geq
0}\|u(t)\|_a\leq c. \ee

 Let us first prove (in paragraph a) below) estimate \eqref{e16}, the 
existence of $u(\infty)$, and
the relation \eqref{e15}, assuming the local existence of the solution
to \eqref{e13}. 

In paragraph b) below the local existence of the solution
to \eqref{e13} is proved.
  
a) If $u(t)$ exists locally, then the function
\be\label{e17}
g(t):=\|\phi\|_{a+\delta}:=\|F(u(t))-h\|_{a+\delta}\ee satisfies the
relation \be\label{e18} g
\dot{g}=(F'(u(t))\dot{u},\phi)_{a+\delta}=-g^2,\ee where equation
\eqref{e13} was used. Since $g\geq 0$, it follows from \eqref{e18}
that \be\label{e19}
 g(t)\leq g(0)e^{-t}, \qquad  g(0)=\|F(u_0)-h\|_{a+\delta}.\ee
From \eqref{e13}, \eqref{e18} and \eqref{e8} one gets:
\be\label{e20} \|\dot{u}\|_a\leq
\frac{1}{c_0}\|\phi\|_{a+\delta}=\frac{g(0)}{c_0}e^{-t}:=re^{-t},\qquad
r:=\frac{\|F(u_0)-h\|_{a+\delta}}{c_0}. \ee
Therefore,
\be\label{e21} \lim_{t\to \infty}\|\dot{u}(t)\|_a=0,\ee and
\be\label{e22} \int_0^\infty\|\dot{u}(t)\|_a dt<\infty. \ee
This inequality  implies 
$$||u(\tau)-u(s)||\leq 
\int_s^{\tau}||\dot{u}(t)dt||<\epsilon, \qquad \tau>s>s(\epsilon),$$
where $\epsilon>0$ is an arbitrary small fixed number, and
$s(\epsilon)$ is a sufficiently large number. Thus, the limit
$V:=\lim_{t\to \infty} u(t):=u(\infty)$ exists by  
the Cauchy criterion, and \eqref{e14} holds.
Assumptions \eqref{e7} and \eqref{e8} and
relations \eqref{e13}, \eqref{e14}, and \eqref{e21} imply \eqref{e15}.

Integrating
inequality \eqref{e20} yields \be\label{e23} \|u(t)-u_0\|_a\leq r,
\ee and \be\label{e24} \|u(t)-u(\infty)\|_a\leq re^{-t}. \ee
Inequality \eqref{e23} implies \eqref{e16}.

b) Let us now prove the local existence of the solution to
\eqref{e13}.

We prove that the operator in
\eqref{e13}
$A^{-1}(u)[F(u)-h]$ is locally Lipschitz in $H_a$.
This implies the local existence of the solution to
\eqref{e13}.
 
One has
\be\label{e25}\begin{split}
&\|A^{-1}(u)(F(u)-h)-A^{-1}(v)(F(v)-h)\|_a\leq
\|[A^{-1}(u)-A^{-1}(v)](F(u)-h)\|_a\\
&+\|A^{-1}(v)(F(u)-F(v))\|_a:=I_1+I_2.\end{split} \ee Write
\be\label{e26} F(u)-F(v)=\int_0^1A(v+t(u-v))(u-v)dt, \ee
and use
assumption \eqref{e9} with $w=v+t(u-v)$ to conclude that
\be\label{e27}I_2\leq c\|u-v\|_a. \ee Write \be\label{e28}
A^{-1}(u)-A^{-1}(v)=A^{-1}(u)[A(v)-A(u)]A^{-1}(v), \ee
and use the
estimate \be\label{e29} \|A^{-1}(v)[F(u)-h]\|_a\leq c,\ee which is a
consequence of assumptions \eqref{e7} and \eqref{e8}. Then use
assumption \eqref{e10} to conclude that 
\be\label{e30} I_1\leq c\|u-v\|_a. \ee 
From \eqref{e25}, \eqref{e27} and 
\eqref{e30} it
follows that the operator $A^{-1}(u)[F(u)-h]$ is 
locally Lipschitz. 

Note that
\be\label{e31} \|u(t)-U\|_a\leq \|u(t)-u_0\|_a+\|u_0-U\|_a\leq
r+\rho, \ee
\be\label{e32} \|F(u(t))-h\|_{a+\delta}\leq
\|F(u_0)-h\|_{a+\delta}\leq
\|F(u_0)-f\|_{a+\delta}+\|f-h\|_{a+\delta}\leq (1+c_0')\rho, 
\ee so, from
\eqref{e20} one gets 
\be\label{e33} r\leq \frac{(1+c_0')\rho}{c_0}. \ee
Choose \be\label{e34} R\geq r+\rho.\ee
 Then the trajectory $u(t)$ stays in the ball $B(U,R)$
for all $t\geq 0$, and, therefore,
 assumptions \eqref{e7}-\eqref{e10} hold in this ball
for all $t\geq 0$.

Condition
 \eqref{e34} and inequality 
\eqref{e33} imply \be\label{e35} \rho\leq \rho_0= \frac{ 
R}{1+c_0^{-1}(1+c_0')}.\ee 
This is
 the "smallness" condition on $\rho.$\\
 Theorem \ref{T:1} is proved. \hfill $\Box$

\section{Example}\label{S:3}

Let
$$F(u)=\int_0^xu^2(s)ds, \qquad x\in [0,1].$$
Then
 $$A(u)q=2\int_0^xu(s)q(s)ds.$$
Let $ f=x$ and $U=1$. Then $F(U)=x$. Choose $a=1$ and $\delta=1$. Denote
by
 $H_a=H_a(0,1)$ the usual Sobolev space. Assume that 
$$h\in B_2(x,\rho):=\{h\ :\ \|h-x\|_2\leq \rho\},$$ 
and  $\rho>0$ is sufficiently small.
One can verify that
 $$A^{-1}(u)\psi=\frac{\psi'(x)}{2u(x)}$$ for any $\psi\in H_1$.

 Let us check conditions \eqref{e7}-\eqref{e11} for this example.

Condition \eqref{e7} holds, because if $u_n\to u$ in $H_1$, then
 $$\int_0^xu_n^2(s)ds\to \int_0^xu^2(s)ds$$ 
in $H_2$.
To verify this, it is sufficient
 to check that
$$\frac{d^2}{dx^2}\int_0^xu_n^2(s)ds\to 2uu',$$
where $\to$ means the convergence in $H:=H_0:=L^2(0,1)$.
In turn, this is verified if one checks that $u_n'u_n\to u'u$ in
$L^2(0,1)$,
 provided that $u'_n\to u'$ in $L^2(0,1)$.

One has
 $$I_n:=\|u_n'u_n-u'u\|_0\leq
 \|(u_n'-u')u_n\|_0+\|u'(u_n-u)\|_0.$$
Since $\|u'_n\|_0\leq c$,
 one concludes that $\|u_n\|_{L^\infty(0,1)}\leq c_1$ and $\lim_{n\to
 \infty}\|u_n-u\|_{L^\infty}=0$. Thus, 
$$\lim_{n\to \infty}I_n=0.$$

Condition \eqref{e8} holds because $\|u\|_{L^\infty(0,1)}\leq
 c\|u\|_1,$ and
$$ \|\int_0^xu(s)q(s)ds\|_2\leq c\|u'q+uq'\|_0\leq
 c(\|q\|_{L^\infty(0,1)}\|u\|_1+\|u\|_{L^\infty(0,1)}\|q\|_1),$$
so
$$\|\int_0^xu(s)q(s)ds\|_2\leq
 c_0'\|u\|_1\|q\|_1,$$
and
$$\|\int_0^xuqds\|_2\geq \|uq\|_1\geq c_0\|q\|_1,$$ provided
that $u\in B_1(1,\rho)$ and $\rho>0$ is sufficiently small.

Condition \eqref{e9} holds because
$$\|A^{-1}(v)A(w)q\|_1=\|\frac{1}{v(x)}w(x)q\|_1\leq c\|q\|_1,$$
provided that $u,w\in B_1(1,\rho)$ and $\rho>0$ is sufficiently
small.

Condition \eqref{e10} holds because
$$\|A^{-1}(u)\int_0^x(u-v)qds\|_1=\|\frac{u-v}{2u}q\|_1\leq
c\|u-v\|_1\|q\|_1, $$
provided that $u,v\in B_1(1,\rho)$ and $\rho>0$
is sufficiently small.

By Theorem \ref{T:1} the equation
$$F(u):=\int_0^x u^2(s)ds=h,$$ where $\|h-x\|_2\leq \rho$ and $\rho>0$
is sufficiently small, has a solution $V$, 
$$F(V)=h.$$ This solution
can be obtained as $u(\infty)$, where $u(t)$ solves problem
\eqref{e13} and conditions \eqref{e11} and \eqref{e35} hold.

\end{document}